\documentclass{amsart}
\usepackage{amscd,amsthm,amssymb,amsfonts,amsmath,euscript}

\theoremstyle{plain}
\newtheorem{thm}{Theorem}[section]
\newtheorem{lemma}[thm]{Lemma}
\newtheorem{prop}[thm]{Proposition}

\theoremstyle{definition}
\newtheorem{defn}[thm]{Definition}

\theoremstyle{remark}
\newtheorem{remark}[thm]{Remark}

\newcommand{\nc}{\newcommand}

\def\makeop#1{\expandafter\def\csname#1\endcsname
  {\mathop{\rm #1}\nolimits}\ignorespaces}
\makeop{Hom}   \makeop{End}   \makeop{Aut}   \makeop{Isom}  \makeop{Pic} 
\makeop{Gal}   \makeop{ord}   \makeop{Char}  \makeop{Div}   \makeop{Lie} 
\makeop{PGL}   \makeop{Corr}  \makeop{PSL}   \makeop{sgn}   \makeop{Spf}
\makeop{Spec}  \makeop{Tr}    \makeop{Nr}    \makeop{Fr}    \makeop{disc}
\makeop{Proj}  \makeop{supp}  \makeop{ker}   \makeop{im}    \makeop{dom}
\makeop{coker} \makeop{Stab}  \makeop{SO}    \makeop{SL}    \makeop{SL}
\makeop{Cl}    \makeop{cond}  \makeop{Br}    \makeop{inv}   \makeop{rank}
\makeop{id}    \makeop{Fil}   \makeop{Frac}  \makeop{GL}    \makeop{SU}
\makeop{Trd}   \makeop{Sp}    \makeop{Tr}    \makeop{Trd}   \makeop{diag}
\makeop{Res}   \makeop{ind}   \makeop{depth} \makeop{Tr}    \makeop{st}
\makeop{Ad}    \makeop{Int}   \makeop{tr}    \makeop{Sym}   \makeop{can}
\makeop{length}\makeop{SO}    \makeop{torsion} \makeop{GSp} \makeop{Ker}
\def\makebb#1{\expandafter\def
  \csname bb#1\endcsname{{\mathbb{#1}}}\ignorespaces}
\def\makebf#1{\expandafter\def\csname bf#1\endcsname{{\bf
      #1}}\ignorespaces} 
\def\makegr#1{\expandafter\def
  \csname gr#1\endcsname{{\mathfrak{#1}}}\ignorespaces}
\def\makescr#1{\expandafter\def
  \csname scr#1\endcsname{{\EuScript{#1}}}\ignorespaces}
\def\makecal#1{\expandafter\def\csname cal#1\endcsname{{\mathcal
      #1}}\ignorespaces} 

\def\doLetters#1{#1A #1B #1C #1D #1E #1F #1G #1H #1I #1J #1K #1L #1M
                 #1N #1O #1P #1Q #1R #1S #1T #1U #1V #1W #1X #1Y #1Z}
\def\doletters#1{#1a #1b #1c #1d #1e #1f #1g #1h #1i #1j #1k #1l #1m
                 #1n #1o #1p #1q #1r #1s #1t #1u #1v #1w #1x #1y #1z}
\doLetters\makebb   \doLetters\makecal  \doLetters\makebf
\doLetters\makescr 
\doletters\makebf   \doLetters\makegr   \doletters\makegr
     \def\qed{\qedmark\medbreak}%
\def\qedmark{{\enspace\vrule height 6pt width 5pt depth 1.5pt}}%

\normalsize

\makeop{Bl}

\def\Fpbar{\overline{\bbF}_p}
\def\Fp{{\bbF}_p}

\newcommand{\Z}{\mathbb Z}
\newcommand{\Q}{\mathbb Q}
\newcommand{\R}{\mathbb R}

\newcommand{\D}{\mathbf D}    
\newcommand{\A}{\mathbb A}    
\newcommand{\F}{\mathbb F}

\newcommand{\npr}{\noindent }

\newcommand{\<}{\langle}   
\renewcommand{\>}{\rangle} 

\newcommand{\isoto}{\stackrel{\sim}{\to}}
\nc{\embed}{\hookrightarrow}

\newcommand{\ch}{characteristic }
\newcommand{\ac}{algebraically closed }
\newcommand{\dieu}{Dieudonn\'{e} }

\nc{\ol}{\overline}
\nc{\wt}{\widetilde}
\nc{\opp}{\mathrm{opp}}
\def\ul{\underline}

\makeop{Ram}
\makeop{Rep}

\begin{document}
\renewcommand{\thefootnote}{\fnsymbol{footnote}}
\setcounter{footnote}{-1}
\numberwithin{equation}{section}


\title{Simple mass formulas on Shimura varieties of PEL-type}
\author{Chia-Fu Yu}
\address{
Institute of Mathematics \\
Academia Sinica \\
128 Academia Rd.~Sec.~2, Nankang\\ 
Taipei, Taiwan and NCTS (Taipei Office)}
\email{chiafu@math.sinica.edu.tw}
\address{
Max-Planck-Institut f\"ur Mathematik \\
Vivatsgasse 7 \\
Bonn, 53111\\ 
Germany}

\date{June 22, 2007. The research is partially supported by 95-2115-M-001-004.}

\begin{abstract}
  We give a unified formulation of a mass for arbitrary abelian
  varieties with PEL-structures and show that it equals a weighted class
  number of a reductive $\Q$-group $G$ relative to an open compact
  subgroup $U$ of $G(\A_f)$, or simply called an {\it arithmetic
  mass}. We classify the special objects for which our formulation
  remains valid over \ac fields. As a result, we show that the set of 
  basic points in a mod $p$ moduli space of PEL-type with a local condition
  (and a mild condition subject to the Hasse principle) can be
  expressed as a double coset space and its mass equals an
  arithmetic mass. The moduli space does not need to have good
  reduction at $p$. This generalizes a well-known
  result for superspecial abelian varieties. 
\end{abstract} 

\maketitle


\def\mass{\mathrm{mass}}
\def\qisom{{\rm Q\text{-}isom}}
\newcommand{\double}[2]{#1(\Q)\backslash #1(\A_f)/#2}
\newcommand{\doublecoset}[1]{\double{#1}{#1(\hat \Z)}}
\def\GAut{{\rm GAut}}
\def\char{{\rm char\,}}

\section{Introduction}
\label{sec:01}

The idea of using elliptic curves to study class numbers 
might go back to Kronecker or even to Gauss. The celebrated 
Eichler-Deuring mass formula says
\begin{equation}
  \label{eq:11}
  \sum_{E\in \Lambda_p} \frac{1}{\#\Aut(E)} =\frac{p-1}{24},
\end{equation}
where $\Lambda_p$ is the set of isomorphism classes of supersingular
elliptic curves over $\ol \F_p$. If we let ${\bfB}:=B_{p,\infty}$ be
the definite quaternion algebra over $\Q$ of discriminant $p$ and
$O_\bfB$ a maximal order, then 
$\# \Lambda_p$ equals the class number of ${\bfB}$. 
Using the adelic language, let $G'$ be
the group scheme over $\Z$ attached to the multiplicative 
group $O_\bfB^\times$; then one has the following natural bijection:
\begin{equation}
  \label{eq:12}
  \Lambda_p\simeq \doublecoset{G'},
\end{equation}
where $\hat \Z=\prod_{p} \Z_p$ and $\A_f=\hat \Z\otimes_{\Z} \Q$ is
the finite ad\`ele ring of $\Q$.
 
We write $\mass(\Lambda_p)$ for the left hand side of (\ref{eq:11}), and
call it the (geometric) mass of $\Lambda_p$. One shows that
$\mass(\Lambda_p)$ equals an arithmetically 
defined mass for the group $G'$ relative to the open compact subgroup
$G'(\hat\Z)$. The definition is as follows. 
For an $\R$-anisotropic reductive $\Q$-group
$G$, and an open compact subgroup $U$ of $G(\A_f)$, the
(arithmetic) mass for $G$ relative to $U$ is defined to be
\begin{equation}
  \label{eq:13}
  \mass(G,U):=\sum_{c} \frac{1}{\# \Gamma_c}, 
\end{equation}
where $c$ runs through a complete set of representatives for the double
coset space $G(\Q)\backslash G(\A_f)/U$, and $\Gamma_c:=G(\Q)\cap
cUc^{-1}$. 

The analogous result of the Eichler-Deuring mass formula for the Siegel
moduli spaces is generalized in Ekedahl 
\cite{ekedahl:ss}. Let $\Lambda_{g,p}$
be the set of the isomorphism classes of $g$-dimensional 
{\it superspecial} principally polarized abelian varieties over 
$\ol \F_p$.  An abelian variety over an \ac field $k\supset \Fp$ 
is called {\it supersingular} if it is
isogenous to a product of supersingular elliptic curves. It is called
{\it superspecial} if it is isomorphic to a product of supersingular
elliptic curves. It is well-known that  $\Lambda_{g,p}$ is a finite
set. Let $G'_g$ be the
group scheme over $\Z$ obtained from $M_g(O_\bfB)$ with equations
given by $g^*g=1$, where $g \mapsto g^*$ is the standard involution.
Similar to the bijection (\ref{eq:12}), one has (see
\cite{ekedahl:ss}, also see \cite[Theorem~2.10,
p.~144]{ibukiyama-katsura-oort})  
\begin{equation}
  \label{eq:14}
    \Lambda_{g,p}\simeq \doublecoset{G_g'}\quad\text{and}\quad
    \mass(\Lambda_{g,p})=\mass(G'_g,G'_g(\hat\Z)), 
\end{equation}
where \[ \mass(\Lambda_{g,p}):=\sum_{(A,\lambda)\in\Lambda_{g,p}}
  \frac{1}{\#\Aut(A,\lambda)}.\]
Applying a formula for $\mass(G'_g,G'_g(\hat\Z))$ in 
Hashimoto-Ibukiyama~\cite[Proposition 9, 
p.~568]{hashimoto-ibukiyama:classnumber} to
the second formula in (\ref{eq:14}), 
Ekedahl obtained the geometric mass formula
\begin{equation}
  \label{eq:15}
  \sum_{(A,\lambda)\in\Lambda_{g,p}}
  \frac{1}{\#\Aut(A,\lambda)}=
  \frac{(-1)^{g(g+1)/2}}{2^g} \left \{ \prod_{k=1}^g \zeta(1-2k)
  \right \}\cdot \prod_{k=1}^{g}\left\{(p^k+(-1)^k\right \}. 
\end{equation}

To generalize geometric mass formulas to Shimura varieties of
PEL-type, the first step is to look for a set $\Lambda$ of certain
special points, such as $\Lambda_{g,p}$ above, in
the reduction $\calM\otimes \Fpbar$ modulo $p$ of a moduli space 
$\calM$ of PEL-type so that
$\Lambda$ can be parameterized by a double coset space of a reductive
group, such as $G'_g$ above. Naturally 
one would consider good reduction cases (as little is known in bad
reduction cases) and consider the set of basic abelian
varieties which are 
``minimal'' in the sense of Oort \cite{oort:minimal}. 
Basic abelian varieties, defined in Kottwitz
\cite{kottwitz:isocrystals}, are those which
land in a minimal Newton polygon stratum of the moduli space 
$\calM\otimes \Fpbar$. 
The minimal basic points in $\calM\otimes \Fpbar$ play
the same role as superspecial points in the Siegel modular varieties;
particularly they form a finite set. Note that the basic locus has
positive dimension in general; thus there is no natural definition
for the mass of the basic locus.
  
It turns out that the correspondence in (\ref{eq:14}) can be more
flexible, even in the Siegel modular varieties. We do not need to assume
that polarizations are principal nor
that supersingular abelian varieties are ``minimal'', that is,
superspecial. 
All we need to do is imposing a type of 
isomorphism class on the associated polarized
$p$-divisible groups. This is proved in \cite{yu:thesis} for
the Hilbert-Blumenthal varieties and in \cite{yu:mass_hb} for the
Hilbert-Siegel modular varieties. It is conceivable that the
same correspondence should hold for basic points in a PEL-type Shimura
variety. 
 

In this paper we give a unified formulation of the (geometric) mass
$\mass(\Lambda)$ for arbitrary abelian varieties with additional
structures over arbitrary (finitely generated) fields. We also show 
that it equals the arithmetic mass defined by a pair $(G,U)$ which can
be explicitly described; see Section~\ref{sec:02} for precise statements. 
The description, though being surprisingly simple, is by no means
obvious. It relies on the deep results of Zarhin, Faltings and de Jong
on the endomorphisms of abelian 
varieties, Tate modules, and $p$-divisible groups; see
\cite{zarhin:end}, \cite{faltings:end} (cf.~\cite{faltings:rational}) and
\cite[Theorem 2.6]{dejong:homo}. We call the formula established in
Theorem~\ref{22} the {\it simple mass
  formula}. The simple mass formula connects masses between the
geometric side and the arithmetic side; but
it provides no clue of computing either side explicitly. It serves as a
comparison theorem between masses arising from different natures in a 
very general setting.   
For example, one can use it to prove a geometric mass formula starting 
from a known arithmetic mass formula and vice versa, or to verify 
an arithmetic mass formula by geometry and vice versa. 
The geometric mass formula~(\ref{eq:15}) is an example. Another
example is given in \cite{yu:mass_hb} where the mass for 
superspecial abelian varieties with real multiplication is
calculated. 
It is  worth to note that a geometric mass then shares good
properties as an arithmetic mass does. For example, it has a simple
relation between different levels and it can be calculated through
local volume computation. 

In some sense the simple mass formula says arithmetic
properties of abelian varieties. For generic abelian varieties, the
simple mass formula provides no information; thus those are not
interesting cases. In the extreme cases such as CM abelian varieties or
supersingular abelian varieties, the simple mass formula provides
non-trivial information and becomes interesting.
In Section 3 we study a class of special abelian
varieties in question which we call {\it of arithmetic type}. Those
include CM abelian varieties and supersingular abelian
varieties. Strictly speaking, this notion is for an abelian variety
$\ul A=(A,\lambda,\iota)$ with additional structures, not for the abelian
variety $A$ itself. Furthermore, it is a ``geometric'' property that is
independent of ground fields over which the object is defined.
The (equivalent) definitions are given in Def.~\ref{32} and Def. \ref{310}. 
For these abelian varieties the hidden Galois structure required in
the mass formula is superfluous; thus the description in 
Theorem~\ref{22} can be extended in a 
geometric setting. Namely, the ground field can be an \ac field; see
Theorem~\ref{314}. 
This explains why a good formulation of the mass for {\it
  supersingular} abelian varieties in \ch $p$ (more generally,
basic abelian varieties with additional structures) or CM abelian varieties
in \ch zero is possible. Some detail 
discussion is included for motivating the
definitions and clarifying the notion of abelian varieties of 
arithmetic type.   
We remark that the description for CM abelian varieties in
\ch zero is well-known. This has been playing an important role on 
explicit reciprocity laws in class field theory, known as the main 
theorem of complex multiplication. 
It is interesting to describe the Galois action on  
the class space in terms of Hecke translation in our setting. 
This might lead an interesting explicit reciprocity law.
                                
In the last section we classify abelian varieties $\ul A$ with additional
structures which are of arithmetic type. In the case of \ch zero, 
the possibility occurs only when every simple factor of the
semi-simple algebra $(B,*)$ is of second kind. Under this condition,  
every object $\ul A$ of arithmetic type is
essentially a self product of a simple 
CM abelian variety (Proposition~\ref{43}). 
In the case of \ch $p$, we show that an object $\ul A$ is 
of arithmetic type if and only if it is basic in the sense of Kottwitz
\cite{kottwitz:isocrystals} (Proposition~\ref{46}). Applying
the results in Sections~\ref{sec:02} and \ref{sec:03}, we complete the
first step of a geometric mass formula for basic points in
Shimura varieties of PEL-type; see Theorem~\ref{main}. \\


\npr {\it Acknowledgments.}
The author thanks J.-K. Yu for helpful discussions in an early stage of
this work. Main part of the manuscript is prepared during the author's
stay at MPIM in Bonn in the fall of 2005. 
He acknowledges the Institute for kind hospitality and 
excellent working environment. \\

\section{Simple mass formulas}
\label{sec:02}

\subsection{}
\label{21}

Let $B$ be a finite-dimensional semi-simple algebra over $\Q$ with a
positive involution $*$, and $O_B$ be any order of $B$ stable under
$*$. 

A {\it polarized abelian $O_B$-variety} is a
triple $\ul A=(A,\lambda,\iota)$, where $(A,\lambda)$ is a polarized abelian
variety and $\iota:O_B\to \End(A)$ is a ring monomorphism such that
$\lambda \iota(a^*)=\iota(a)^t \lambda$ for all $a\in O_B$. 
For any $\ul A$ and any prime $\ell$ (not necessarily invertible in
the ground field), we write $\ul A(\ell)$ for the associated
$\ell$-divisible group with additional structures $(A[\ell^\infty],
\lambda_\ell, \iota_\ell)$, where $\lambda_\ell$ is the induced
quasi-polarization from $A[\ell^\infty]$ to
$A^t[\ell^\infty]=A[\ell^\infty]^t$ (the Serre dual), 
and $\iota_\ell:O_B\otimes\Z_\ell\to \End(A[\ell^\infty])$ 
the induced ring monomorphism. 

For any two $\ul A_1$ and $\ul A_2$ over a field $k$,
denote by 
\begin{itemize}
\item $\qisom_k(\ul A_1, \ul A_2)$ (resp.~$\Isom_k(\ul A_1, \ul
A_2)$) the set of $O_B$-linear 
quasi-isogenies (resp.~isomorphisms) $\varphi:A_1\to A_2$ over $k$
such that $\varphi^*\lambda_2=\lambda_1$, and   
\item $\qisom_k(\ul A_1(\ell), \ul A_2(\ell))$ (resp. $\Isom_k(\ul
A_1(\ell), \ul A_2(\ell))$) the set of $O_B\otimes \Z_\ell$-linear 
quasi-isogenies (resp. isomorphisms) $\varphi:A_1[\ell^\infty]\to
A_2[\ell^\infty]$ such that $\varphi^*\lambda_2=\lambda_1$.
\end{itemize}

Let $k$ be a field of finite type over
its prime field, and    
let $x:=\ul A_0=(A_0,\lambda_0,\iota_0)$ be a fixed polarized abelian
$O_B$-variety over $k$.
Denote by $\Lambda_x(k)$ the set of isomorphisms classes of polarized
abelian $O_B$-varieties $\ul A$ over $k$ such that 
\begin{itemize}
\item [(i)] ($I_\ell$): $\Isom_k(\ul
A_0(\ell),\ul A(\ell))\neq \emptyset$ for all $\ell$, and
\item [(ii)] $(Q)$: $\qisom_k(\ul A_0,\ul A)\neq \emptyset$.
\end{itemize}
Let $G_{x}$ be the
automorphism group scheme over $\Z$ associated to $\ul A_0$; for any
commutative ring $R$, its group of $R$-valued points is
\begin{equation}
  \label{eq:21}
  G_{x}(R)=\{g\in \End_{O_B}(A_{0/k})\otimes R\, | \, g' g=1\},
\end{equation}
where $g\mapsto g'$ is the Rosati involution induced by
$\lambda_0$. By definition, \[ G_{x}(\Q)=\qisom_k(\ul A_0, \ul
A_0). \] 
For the group of $\Q_\ell$-valued points of $G_x$, we have

\begin{thm}[Zarhin, Faltings, de Jong]\label{215}
  Notations as above,
  one has the natural isomorphisms
\begin{equation}
  \label{eq:22}
  G_{x}(\Z_\ell)=\Isom_k(\ul A_0(\ell),\ul A_0(\ell))\quad 
\text{and}\quad G_{x}(\Q_\ell)=\qisom_k(\ul A_0(\ell),\ul A_0(\ell))
\end{equation}
for all $\ell$.
\end{thm}
\begin{proof}
  This follows immediately from the theorem of Tate type:
\[ \Hom_{k}(A_1,A_2)\otimes \Z_\ell\simeq
\Hom_k(A_1[\ell^\infty],A_2[\ell^\infty]), \]
for all $\ell$. 
This is due to Zarhin in the case where $\char k=p\neq \ell$, due to
Faltings in the case where $\char k=0$, and
due to de Jong in the case where $\char k=p=\ell$.
See \cite{zarhin:end}, \cite{faltings:end}
(cf.~\cite{faltings:rational}) and  
\cite[Theorem 2.6]{dejong:homo} for more details. \qed 
\end{proof}

\begin{thm}\label{22}

{\rm (1)} There is a natural bijection between the following two
pointed sets 
\[ \Lambda_x(k)\simeq  
  G_{x}(\Q)\backslash G_{x}(\A_f)/G_{x}(\hat \Z). \]
Consequently, $\Lambda_x(k)$ is finite.

{\rm (2)} Define \[ \mass[ \Lambda_x(k)]:=\sum_{\ul A\in
  \Lambda_x(k)} \frac{1}{\#\Aut_k(\ul A)}, \]
where $\Aut_k(\ul A):=\Isom_k(\ul A, \ul A)$ defined as before. 
Then one has $\mass[ \Lambda_x(k)]=\mass[G_{x},
G_{x}(\hat \Z)]$. 
\end{thm}
\begin{proof}
  (1) Given an element $\ul A\in \Lambda_x(k)$, consider the natural map
   \begin{equation}
    m(\ul A): \label{eq:23}
    \qisom(\ul A,\ul A_0)\times \prod_\ell \Isom_k(\ul A_0(\ell),\ul
    A(\ell))  
    \to \prod_\ell\,\!\! '\, \qisom_k(\ul A_0(\ell),\ul A_0(\ell))=
    G_{x}(\A_f)  
  \end{equation}
which sends $(\phi,(\alpha_\ell)_\ell)$ to 
$(\phi \alpha_\ell)_\ell$. If $c$ is an element
in the image of $m(\ul A)$, then the image  $c(\ul A)$ is
the double coset $G_{x}(\Q)\, c\, G_{x}(\hat \Z)$. Thus, the map 
$\ul A \mapsto c(\ul A)$
defines a map from $\Lambda_x(k)$ to the double coset space
$\doublecoset{G_{x}}$. 

Let $\ul A, \ul A'\in \Lambda_x(k)$ such that $c(\ul A)=c(\ul
A')$. Write $c(\ul A)=[(\phi \alpha_\ell)_\ell]$ and   $c(\ul
A')=[(\phi' \alpha'_\ell)_\ell]$. Then there exist $b\in G_x(\Q)$ and
$k_\ell\in G_x(\Z_\ell)$ for all $\ell$ such that $b\phi \alpha_\ell
k_\ell=\phi' \alpha'_\ell$. Then 
\[ (b\phi)^{-1} \phi'=\alpha_\ell
k_\ell (\alpha'_\ell)^{-1}\in \qisom_k(\ul A',\ul A)\cap \prod_{\ell}
\Isom_k(\ul A'(\ell),\ul A(\ell))=\Isom_k(\ul A',\ul A).\] 
Thus $\ul A'\simeq \ul A$ and this shows the injectivity of $c$.

Given $[(\phi_\ell)_\ell]$ in $\doublecoset{G_x}$, choose a positive
integer $N$ such that $f_\ell:=N\phi^{-1}_\ell$ (from
$A_0[\ell^\infty]$ to itself) is an isogeny for all
$\ell$. Let $H$ be the product of the kernels of $N\phi^{-1}_\ell$; this
is a finite subgroup scheme over $k$ invariant under the $O_B$-action. 
Take $A:=A_0/H$ and let $\pi:A_0\to A$ be the natural projection; $A$
is defined over $k$ and is equipped with a natural action by $O_B$
so that $\pi$ is $O_B$-linear. Let $\lambda\in \Hom(A,A^t)\otimes \Q$
be the fractional polarization on $A$ such that
$(N^{-1}\pi)^*\lambda=\lambda_0$; $\lambda$ is $O_B$-linear as $\pi$ is
so. As $\pi_\ell$ and $f_\ell$ have the same kernel, there is an
element $\alpha_\ell\in \Isom_{O_B,k}(\ul A_0[\ell^\infty],
\ul A[\ell^\infty])$ such that $\alpha_\ell f_\ell =\pi_\ell$. 
Then $\alpha_\ell= \pi_\ell f_\ell^{-1}=N^{-1} \pi_\ell \phi_\ell$ and
$\alpha_\ell^* \lambda=\lambda_0$. 
This shows $\lambda\in \Hom_{k,O_B}(A,A^t)$ and one obtains 
an object $\ul A\in \Lambda_x(k)$. Put 
$\phi:=(N^{-1}\pi)^{-1}\in \qisom_k(\ul A,\ul A_0)$. One checks 
\[ \phi \alpha_\ell=N \pi_\ell^{-1} \alpha_\ell=N f_\ell^{-1}
=\phi_\ell. \]  
This shows $c(\ul A)=[(\phi_\ell)_\ell]$ and the surjectivity of the
map $c$. 

(2) It suffices to show that if $x'=\ul A\in \Lambda_x(k)$ and
  $c$ any representative for the double coset $c(\ul
A)$, then $\Aut_k(\ul A)\simeq \Gamma_c$. 
Write $G_{x'}$ for the group scheme over $\Z$ associated to $\ul A$
defined as in \ref{21}. 
Choose $\phi\in \qisom_k(\ul A_0,\ul A)$
such that $\phi c_\ell\in \Isom_k(\ul A_0(\ell),
  \ul A(\ell))$ for all $\ell$. 
  Note that $\alpha\in
  \Aut_k(\ul A)$ if and only if $\alpha\in G_{x'}(\Q)$ and
  $\alpha_\ell\in \Aut_k(\ul A(\ell))$ for all $\ell$. 

  The map $\phi$ induces an isomorphism $G_x(\Q)\isoto G_{x'}(\Q)$
  which sends 
  $\beta$ to $\phi\beta\phi^{-1}=:\alpha$. Note that
  $\alpha\in G_{x'}(\hat \Z)$ if and only if $(\phi c)^{-1} \alpha 
  (\phi c)\in G_x(\hat \Z)$. The latter is equivalent to 
  $c^{-1}\beta c\in G_x(\hat \Z)$. Therefore, the above isomorphism 
  gives $ \Gamma_c \simeq \Aut_k(\ul A)$. This completes the
  proof. \qed  
\end{proof}

\subsection{}
\label{23}
Let $N$ be any positive integer and $U_N$ be the kernel of the
reduction map $G_x(\hat \Z)\to G_x(\hat \Z/N\hat \Z)$. Let $\ul A$ be
a polarized abelian $O_B$-variety. By an {\it $(\ul A_0,U_N)$-level
structure on $\ul A$} we mean a non-empty $U_N$-orbit $\bar \eta$ of
isomorphisms $\eta$ in $\prod_\ell \Isom_k(\ul A_0(\ell),\ul
A(\ell))$. The 
existence of such $\bar \eta$ implies that the condition (i) for
$\Lambda_x(k)$ is satisfied. Let $\bar \eta_0$ be the 
$U_N$-orbit of the identity in $\prod_\ell \Isom_k(\ul A_0(\ell),\ul
A_0(\ell))$. We change our notation a bit in the remaining of this
section. We write $\ul A_0$ for $(A_0,\lambda_0,\iota_0, \bar \eta_0)$
and $\ul A$ for $(A,\lambda,\iota, \bar \eta)$ in brief.


For any two $\ul A_1$ and $\ul A_2$ over a field $k$, denote by 
\begin{itemize}
\item $\qisom_k(\ul A_1, \ul A_2)$ and $\qisom_k(\ul A_1(\ell), \ul
A_2(\ell))$ the sets which have the same meaning as in (\ref{21});
\item $\Isom_k(\ul A_1, \ul A_2)$ the set of elements $\varphi$ in 
$\Isom_k((A_1,\lambda_1,\iota_1),(A_2,\lambda_2,\iota_2))$ satisfying
$\varphi_* \bar \eta_1=\bar \eta_2$;  
\item $\Isom_k(\ul
A_1(\ell), \ul A_2(\ell))$ the set of elements $\varphi$ in  
$\Isom_k((A_1,\lambda_1,\iota_1)(\ell),(A_2,\lambda_2,\iota_2)(\ell))$
satisfying $\varphi_* \bar \eta_{1,\ell}=\bar \eta_{2,\ell}$.
\end{itemize}


Let $\Lambda_{x,N}(k)$ denote the set of isomorphism classes of
polarized abelian $O_B$-varieties with an $(\ul A_0,U_N)$-level
structure $(A,\lambda,\iota,\bar \eta)$ over $k$ such that 
$\qisom_k(\ul A_0,\ul A)\neq \emptyset$. 
The same proof of Theorem~\ref{22} without modification gives the
following variant.


\begin{thm}\label{24}
There is a natural bijection 
\[ \Lambda_{x,N}(k)\simeq  
  G_{x}(\Q)\backslash G_{x}(\A_f)/U_N. \]
Furthermore, one has $\mass[ \Lambda_{x.N}(k)]=\mass(G_{x},
U_N)$.
\end{thm}

\begin{lemma}
  If $N\ge 3$, then the group $\Aut_k(\ul A)$ is trivial for any object $\ul
  A=(A,\lambda,\iota,\bar \eta)$ in $\Lambda_{x,N}(k)$. 
\end{lemma}
\begin{proof}
Clearly $\Aut_k(\ul A)\subset \Aut_k(A,\lambda,\bar \eta)$. By a
theorem of Serre, we have $\Aut_k(A,\lambda,\bar
\eta)=\{1\}$. Therefore, $\Aut_k(\ul A)$ is trivial. \qed
\end{proof}

\begin{remark}
  When the ground $k$ has \ch $p$, the usual assumption of $O_B$ for
  good reduction is not required in Theorems~\ref{22}
  and \ref{24}.
\end{remark}

\section{Abelian varieties of arithmetic type}
\label{sec:03}

In this section we define a class of abelian varieties with additional
structures those share the similar properties as supersingular abelian
varieties. We call them of arithmetic type (or of $B$-arithmetic type
to be precise) as those have rich arithmetic properties. 
The naive definition is given in
Def. \ref{32}. After studying basic properties of these abelian
varieties, we give another equivalent but refined definition in
Def. \ref{310}.

\subsection{}
\label{31}
In the following, polarized $O_B$-abelian varieties will be assumed
implicitly to have the properties:
($*$) There is a non-degenerate
$\Q$-valued skew-Hermitian $B$-space $(V,\psi)$ such that
$2 \dim A=\dim_{\Q} V$.  

That is, $V$ is a left faithful finite $B$-module, and 
$\psi:V\times V\to \Q$ is a non-degenerate alternating pairing on $V$
such that $\psi(bx,y)=\psi(x,b^* y)$ for all $b\in B$ and $x,y,\in V$.
In other words, the objects $\ul A$ we consider are those that 
can occur in the Shimura variety attached to the reductive group
defined by a PEL datum $(B,*, V,\psi)$. 
However, the moduli space of PEL-type itself will 
not be specified explicitly here. 


\begin{defn}\label{32} Let $k$ be a field of finite type over its
  prime field. Let $\ul A=(A,\lambda,\iota)$ be a polarized abelian
  $O_B$-variety over $k$ such that $\End(A_{\bar
  k})=\End(A_k)$. 
Let $T_\ell:=T_\ell (A)$ denote the Tate module of 
$A$ ($\ell \neq \char k$). Write
$\GAut_B(T_\ell)$ or $\GAut_B(T_\ell, \<\, ,\>)$ for the group 
\[ \{g\in \End_{B_\ell}(V_\ell);\ g' g\in \Q_\ell^\times\}, \]
where $g\mapsto g'$ is the involution induced by $\lambda$. 
Let $\rho_\ell:\calG_k:=\Gal(k_s/k)\to \GAut_B(T_\ell)$ 
be the associated $\ell$-adic Galois representation. 
We call $\ul A$ over $k$ is {\it of arithmetic type} 
if the image $\rho_\ell( \calG_k)$ is contained 
in the center of $\GAut_B(T_\ell)$ 
for all $\ell\neq \char k$.   
\end{defn}

We write $B_\ell:=B\otimes \Q_\ell$ and $V_\ell=T_\ell\otimes
\Q_\ell$. Let $G_\ell:=\rho_\ell(\calG_k)$, the image of the
$\ell$-adic Galois representation, and let $G^{\rm alg}_\ell$ be the
algebraic envelope of $G_\ell$, the smallest algebraic subgroup $H$ of
$\GL(V_\ell)$ over $\Q_\ell$ such that $H(\Q_\ell)$ contains $G_\ell$.
Write $G^0$ for the neutral component of an algebraic group $G$.


 

 

 

\begin{lemma}\label{33}
  Let $\ul A=(A,\lambda,\iota)$ be a polarized abelian $O_B$-variety
  over $k$ as in Def.~\ref{32}.

{\rm (1)} If $O_B=\Z$, that is, $\ul A=(A,\lambda)$ be a polarized
  abelian variety, then $\ul A$ is of arithmetic type if and only if
$\char k=p>0$ and $A$ is supersingular.

{\rm (2)} If $B$ is a commutative semi-simple algebra and
  $[B:\Q]=2\dim A$, that is, $A$ has complex multiplication by $B$,
  then $\ul A$ is of arithmetic type. 

{\rm (3)} Let $\ul A'$ is another polarized
abelian $O_B$-variety such that there is an $O_B$-linear isogeny
  $\varphi: A\to A'$ which preserves the polarizations. Then $\ul A$
  is of arithmetic type if and only if $\ul A'$ is also of arithmetic
  type. 

{\rm (4)}  Let $B$ be a semi-simple subalgebra of a matrix algebra
     $M_n(K)$, and let $C_B$ be the commutant of $B$ in $M_n(K)$. 
     Then the center $Z(C_B)$ of $C_{B}$ of $B$ is 
     the center $Z(B)$ of $B$.
\end{lemma}
\begin{proof}
  (1) If the image of the
Galois group lies in 
the center of $\GSp_{2g}(\Z_\ell)$, which consists of
the scalar matrices, then one has $\dim \End^0(A)=4g^2$. This is 
possible only when $\char k=p>0$ and $A$ is supersingular.

(2) Since $\End_{O_B}(T_\ell(A))$ is already commutative, the
    statement follows.

(3) This is clear.

(4)  By the bi-commutant theorem, the elements of $M_n(K)$ which commute
  with $C_{B}$ lie in $B$. Therefore, $Z(C_B)$ consists of elements in
  $B$ commuting with $B$. This shows $Z(C_B)=Z(B)$. \qed
\end{proof}

\subsection{}\label{34}


Let $\ul A=(A,\lambda,\iota)$ is a polarized abelian $O_B$-variety
over $k$ as in Def.~\ref{32}. Write the semi-simple algebra $B$ into
simple factors $\oplus_{i=1}^r M_{n_i}(D_i)$, where $D_i$ is a division
algebra over $\Q$ with a positive involution $*_i$. According to this
decomposition the abelian variety $A$ is isogenous to $\prod
A_i^{n_i}$. 
One has ring monomorphisms $\iota_i: D_i\to
\End^0(A_i):=\End(A_i)\otimes \Q$. Write $V_i$
for $T_\ell(A_i)\otimes \Q_\ell$ and one has
\[ V_\ell=\prod_{i=1}^r V_i^{\oplus n_i}, \]
\[ \End_{B}(V_\ell)=\oplus_{i=1}^r \End_{D_i}(V_i). \]
Let $g\mapsto g'$ be the adjoint with respect to the alternating
pairing $\<\,,\>$ on $V_\ell$. Then $\GAut_{B}(V_\ell,\<\,,\>)$
consists of elements $g=(g_i)\in \prod \End_{D_i}(V_i)$ such that
$g'_1 g^{ }_1= g'_2 g^{ }_2=\dots=g'_r g^{ }_r\in \Q_\ell^\times$. 

We have projections 
$p_i: \GAut_{B}(V_\ell,\<\,,\>)\to\GAut_{D_i}(V_i,\<\,,\>_i)$ and
these induce $p_i: Z(\GAut_{B}(V_\ell,\<\,,\>))\to
Z(\GAut_{D_i}(V_i,\<\,,\>_i))$. 
If $\rho_{i,\ell}$ is the $\ell$-adic Galois
representation attached to $\ul A_i$, then one has $p_i\circ
\rho_\ell=\rho_{i,\ell}$. This shows that if $\ul A$ is of arithmetic
type, then so as each factor $\ul A_i$. The converse is also true
as  $Z(\GAut_{B}(V_\ell,\<\,,\>))\hookrightarrow
\prod_i Z(\GAut_{D_i}(V_i,\<\,,\>_i))$. Thus one has

\subsubsection{}\label{341} Notation as above. Then $\ul A$ is of
arithmetic type if and only if each factor $\ul A_i$ is of arithmetic type. \\


\begin{defn}\label{35}
  Let $(B,*)$ be as in \ref{21} and assume that $B$ is a division
  algebra. A polarized abelian $O_B$-variety is said to be of type
  (D \& 0-dim) if $B$ is of type (III) in the Albert classification and 
  $2 \dim A =[B:\Q]$.   
\end{defn}

Recall \cite[Section 23]{mumford:av} that $B$ is of type (III) if $B$
is a totally definite quaternion algebra over a totally real number
field $F$ and $*$ is the main involution, the unique positive involution.


\begin{lemma}\label{36} Assume that $B$ is a division algebra. Let
  $\ul A$ be a polarized abelian $O_B$-variety of type
  {\rm (D \& 0-dim)}. Then 

{\rm (1)} $T_\ell$ is a free $O_F\otimes \Z_\ell$-module of rank
  $4$.

{\rm (2)}  $V_\ell$ is a free $B_\ell$-module of
  rank $1$.

{\rm (3)}  $\GAut_{B_\ell}(V_\ell,\<\,,\>)$ is an extension of a normal
  commutative subgroup by a finite 2-torsion group.

{\rm (4)} The center $Z(\GAut_{B_\ell}(V_\ell,\<\,,\>))$ consists
  of elements $a$ in $F_\ell:=F\otimes \Q_\ell$ with $a^2\in
  \Q_\ell^\times$. 
\end{lemma}
\begin{proof}
  The statement (1) follows from the fact that $\Tr(a;V_\ell/
  \Q_\ell)=4 \Tr_{F/\Q}(a)$ for all $a\in O_F$. The statement (2)
  follows from (1). To show the statement (3), we regard
  $G:=\GAut_{B}(V_\ell,\<\,,\>)$ as an algebraic group over $\Q_\ell$
  and show that its neutral component $G^0$ is a torus. 

  Let $V_\ell=B_\ell$ as a left $B_\ell$-module. 
  Let $(\,,):B_\ell\times B_\ell\to B_\ell$ be the
  lifting of $\<\,,\>$. One has $\<x,y\>=\Trd_{B_\ell/\Q_\ell}(x\alpha
  y^*)$, where $\alpha=(1,1)$ and $\alpha^*=-\alpha$. Any element in
  $\End_{B_\ell}(V_\ell)$ is a right translation $\rho_g$ for some $g\in
  B_\ell$. The condition $\<xg,yg\>=c(g)\<x,y\>$ gives
  $\Trd_{B_\ell/\Q_\ell}(xg\alpha g^* y^*)=\Trd_{B_\ell/\Q_\ell}(x\,
  c(g)\alpha y^*)$. Therefore, the group $G$ is the subgroup of
  $B_\ell^{{\rm opp},\times}$ defined by the relation $g \alpha
  g^*=c(g) \alpha $ for some $c(g)\in \Q_\ell^\times$. 
  Choose the isomorphism $B_\ell^{\rm opp}\simeq B_\ell$
  which sends $g\mapsto g^{-1}$; the group $G$ is identified with the
  subgroup of $B_\ell^\times$ defined by the same relation.  

  For each $\sigma\in \Sigma:=\Hom(F_\ell,\ol \Q_\ell)$, put
  $B_\sigma=B_\ell\otimes_{F_\ell,\sigma} \ol \Q_\ell\simeq
  M_2(\ol \Q_\ell)$. Let $j=
  \begin{pmatrix} 0 & -1 \\ 1 & 0 \end{pmatrix}$ 
  and $g\in B_\sigma$, one computes 
\[ j g^* j^{-1}=\begin{pmatrix} 0 & -1 \\ 1 & 0 \end{pmatrix}
\begin{pmatrix} d & -b \\ -c & a \end{pmatrix}
\begin{pmatrix} 0 & 1 \\ -1 & 0 \end{pmatrix}=
\begin{pmatrix} a & c \\ b & d \end{pmatrix}=g^t. \]
Write $\alpha=\beta j$, then $\beta^t=j \beta^* j^{-1}=-\alpha^*
j^{-1}=\beta$ and the relation defining $G$ becomes $g \beta
g^t=c(g) \beta$ for some $c(g)$. We have proved
\[ G_{\ol \Q_\ell}\simeq \{(g_\sigma)\in \GL_2^\Sigma\ ;\  
g^{}_\sigma g_\sigma^t=c\ \text{ for some } c\in \ol \Q^\times_\ell
\text{\,(independent of $\sigma$)\,},
\forall\, \sigma\in \Sigma\},\ \text{ and} \]
\[ G^0_{\ol \Q_\ell} \simeq \left \{
  \begin{pmatrix} a_i & b_i \\ -b_i & a_i \end{pmatrix}
   \in \GL_2^d\ ;\ a_i^2+b_i^2 =c\ \text{ for some } 
    c\in \ol \Q^\times_\ell, 
   \forall\, 1\le i\le d=[F:\Q] \right \}. \]
This shows that $G^0$ is a torus.

(4) This follows directly from the computation in (3). \qed
\end{proof}

\begin{lemma}\label{37}
For any polarized abelian $O_B$-variety $\ul A$, 
the center
$Z(\GAut_{B}(V_\ell,\<\,,\>))$ consists of elements $a\in Z(B)\otimes
\Q_\ell$ such that $a^* a\in \Q_\ell^\times$, where $\Q_\ell$ is
included in $Z(B)\otimes \Q_\ell$ diagonally.     
\end{lemma}

\begin{proof}
By the argument (\ref{34}) one reduces to the case where $B$ is a
division algebra. 
The case of type (D \& 0-dim) has been treated in Lemma~\ref{36}. 
Now suppose that $\ul A$ is not of type (D \& 0-dim).

For an algebra $E$ and a subset $G$ of $E$, we write 
\[ Z(E,G):=\{x\in E;\ gx=xg\
\forall\,g\in G\}.\]
Let $E:=\End_B(V_\ell)$ and $G$ be the 
algebraic group over $\Q_\ell$ defined by $\ul V_\ell$; we have
particularly $G(\Q_\ell)=\GAut_{B}(V_\ell,\<\,,\>)$. It then suffices
to show that $Z(B)\otimes \Q_\ell=Z(E,G(\Q_\ell))$. We have
$Z(B)\otimes \Q_\ell\subset Z(E,G(\Q_\ell))$. So it suffices to show
$\dim_{\Q_\ell} Z(B)\otimes \Q_\ell= \dim_{\Q_\ell}
Z(E,G(\Q_\ell))$. 
Since $\dim_{\Q_\ell} Z(B)\otimes
  \Q_\ell=\dim_{\ol \Q_\ell} Z(B)\otimes \ol \Q_\ell$ and 
\[ \dim_{\ol \Q_\ell} 
  Z(E\otimes \ol \Q_{\ell},G^0(\ol \Q_\ell))= \dim_{\Q_\ell}
  Z(E,G^0(\Q_\ell))\ge \dim_{\Q_\ell}  
  Z(E,G(\Q_\ell))\] 
(as $G^0(\Q_\ell)$
is Zariski dense in $G^0$), it suffices to show that 
$\dim_{\ol \Q_\ell} Z(B)\otimes \ol \Q_\ell=\dim_{\ol \Q_\ell} 
  Z(E\otimes \ol \Q_{\ell},G^0(\ol \Q_\ell))$. Decomposing into simple
factors, we have three cases:

\begin{enumerate}
  \item[(a)] $E=M_n(\ol \Q_\ell)\times M_n(\ol \Q_\ell)$,
  $*:(A,B)\mapsto (B^t,A^t)$ and $G=GU_n$.
  \item[(b)] $E=M_{2n}(\ol \Q_\ell)$, $*$ is the standard symplectic
  involution, and $G=GSp_{2n}$. 
  \item[(c)] $E=M_{2n}(\ol \Q_\ell)$, $*:A\mapsto A^t$ and $G=GO_{2n}$
  ($n\ge 2$).   
\end{enumerate}
Then we have the cases (a)~$Z(E,G^0)=\{(aI_n, bI_n);a,b\in \ol \Q_\ell\}$;
(b)~$Z(E,G^0)=\{aI_{2n};a\in \ol \Q_\ell\}$; (c)~$Z(E,G^0)=\{aI_{2n};a\in
  \ol \Q_\ell\}$. From this one  sees that 
  $\dim_{\ol \Q_\ell} Z(B)\otimes \ol \Q_\ell=
  \dim_{\ol \Q_\ell} Z(E\otimes \ol \Q_{\ell},G^0(\ol \Q_\ell))$. 
   This finishes the proof. \qed
\end{proof}

\begin{lemma}\label{38}
  Let $\ul A=(A,\lambda,\iota)$ be a polarized abelian $O_B$-variety
  and $\lambda'$ be another $O_B$-linear polarization. Then
  $(A,\lambda',\iota)$ is of arithmetic type if and only if $\ul A$ is
  of arithmetic type. 
\end{lemma}
\begin{proof}
  By Lemma~\ref{37}, the center of $\GAut_{B}(V_\ell,\<\,,\>)$ is
  independent of the choice of polarizations. 
  Therefore, the assertion is proved. \qed 
\end{proof}


\begin{lemma}\label{39}
  Let $k_0$ be a field of finite type over its prime field and $k$ be 
  a field finitely generated over $k_0$. 
  Let $\ul A$ be a polarized
  abelian $O_B$-variety over $k_0$. If $\ul A$ is of arithmetic type
  then so as $\ul A_k:=\ul A\otimes_{k_0} k$. Conversely, if $\ul A_k$ is of
  arithmetic type, then so as $\ul A_{k'}$ for a finite extension $k'$
  of $k_0$. 
\end{lemma}
\begin{proof}
  Let $k_s$ be a separable closure of $k$ and $k_{0,s}$ the algebraic
  closure of $k_0$ in $k_s$. Let $k_1$ be the algebraic closure of
  $k_0$ in $k$.
  The restriction gives a surjective homomorphism $r:\calG_{k}\to
  \calG_{k_1}$ of Galois groups. We also have the Galois equivariant
  isomorphism $s:A[\ell^n](k_{0,s})\simeq A[\ell^n](k_{s})$ in the
  sense that $r(\sigma)x=\sigma(s(x))$ for $x\in A[\ell^n](k_{0,s})$
  and $\sigma\in \calG_k$. This gives rise to the commutative diagram
\[ \begin{CD}
  \calG_{k}  @>\rho_{A_k}>> \Aut (T_\ell(A_{k_s})) \\
  @V r VV  @V\simeq VV  \\
  \calG_{k_1}  @>\rho_{A_{k_0}}>> \Aut (T_\ell(A_{k_{0,s}})), \\
\end{CD} \] 
and one has $\rho_{A_k}(\calG_{k})=\rho_{A_{k_0}}(\calG_{k_1})$. It
follows that $\ul A_k$ is of arithmetic if and only if $\ul 
A_{k_1}$ is so. From this the first and second statement follow. \qed
\end{proof}

Lemma~\ref{33} (3) and Lemmas~\ref{38} and \ref{39} suggest the
following alternative definition, which is more refined.

\begin{defn}\label{310}
  Let $(B,*)$ be as in \ref{21} and $(A,\iota)$ be an abelian
  $B$-variety up to isogeny over $k$ as in Def.~\ref{32}. 
  The pair $(A,\iota)$ is said to be {\it of
  $B$-arithmetic type} or simply {\it of arithmetic type} if there is a
  finite extension $k'/k$ such that $\rho_\ell(\calG_{k'})$ lies in the
  center of $\GAut_{B}(V_\ell)$ for some $B$-linear polarization
  $\lambda$ and for one $\ell\neq \char k$. An abelian $B$-variety up
  to isogeny is said to be {\it of arithmetic type} if it is so over a
  field of finite type over its prime field.  
\end{defn}

For any abelian $B$-variety $(A,\iota)$ up to isogeny, there always exists a
$B$-linear polarization on $(A,\iota)$ (see \cite[Section 9]{kottwitz:jams92}).
We now show the $\ell$-independence of $G_\ell^{\rm alg}$ for abelian
varieties of arithmetic type. 

\begin{lemma}\label{311}
  Let $\ul A=(A,\iota)$ be an abelian $B$-variety up to
  isogeny. If $\ul A$ is of arithmetic type, 
  then $A$ is of CM type.  
\end{lemma}
\begin{proof}
  Since $\ul A$ is of arithmetic type, $G_\ell$ is
  commutative. Denote by $\Q_\ell[G_\ell]$ the subalgebra of $\End(V_\ell)$
  generated by $G_\ell$. By the semi-simplicity of Tate modules
  due to Faltings and Zarhin \cite{faltings:end, zarhin:end},
  $\Q_\ell[G_\ell]$ is a (commutative) 
  semi-simple subalgebra. Let $L$ be a maximal semi-simple commutative 
  subalgebra in $\End^0(A)$, then $L\otimes \Q_\ell$ is a maximal
  commutative semi-simple algebra in $\End^0(A)\otimes \Q_\ell$. By
  the theorem of Faltings and Zarhin on Tate's conjecture loc.~cit., 
  we have $\End^0(A)\otimes \Q_\ell=\End_{\Q_\ell[G_\ell]}(V_\ell)$. Hence
  $L\otimes \Q_\ell$ becomes a maximal semi-simple commutative
  subalgebra in $\End_{\Q_\ell[G_\ell]}(V_\ell)$. Since
  $\Q_\ell[G_\ell]$ is commutative and semi-simple, any maximal
  semi-simple commutative subalgebra in $\End_{\Q_\ell[G_\ell]}(V_\ell)$
  has degree $2g$ over $\Q_\ell$. This shows $[L:\Q]=2g$ and the proof
  is complete. \qed    
\end{proof}

\begin{prop}\label{312}
  Let $(A,\iota)$ be an abelian $B$-variety of arithmetic type over
  $k$ as in Def.~\ref{310}.
  Then 
  the algebraic envelope $G^{\rm alg}_\ell$ is independent of 
  $\ell$ for all $\ell \neq \char k$. That is, there is a
  $\Q$-subgroup $G$ of $\GL_{2g}$ such that
  $G\otimes \Q_\ell \simeq G^{\rm alg}_\ell$ for all 
  $\ell \neq \char k$.  
\end{prop}
\begin{proof}
  By Lemma~\ref{311}, $A$ is of CM-type. The semi-simple part of
  $G_\ell^{\rm alg}$ is trivial. By Bogomolov's theorem, 
  the neutral component $(G^{\rm alg}_\ell)^0$ is independent of $\ell$ (see
  \cite[2.2.5]{serre:1984-5}, 
  also see the remark in 2.3 of loc.~cit.~for the function field
  case). 
  By a theorem of Serre \cite{serre:1984-5} that the component group 
  $G^{\rm alg}_\ell/(G^{\rm alg}_\ell)^0$ is independent of $\ell$,
  one shows that $G^{\rm alg}_\ell$ is independent of $\ell$. \qed  
\end{proof}

\begin{remark} \label{313}
In \ref{33} -- Proposition~\ref{312} we have shown that $\ul A$ is of
arithmetic 
type in the sense of Def.~\ref{32} if and only if 
its underlying abelian $O_B$-variety $(A,\iota)$ is of
arithmetic type in the sense of Def.~\ref{310}. According to Definition
\ref{310}, the property for an object $\ul A$ being of arithmetic type is
independent of the ground field for which it is defined.
\end{remark}

\begin{thm}\label{314} 
Let $k$ be an \ac field and let $x=\ul A$ be a polarized abelian
$O_B$-variety over $k$.  Suppose that $\ul A$ is of arithmetic type.

  {\rm (1)} $G_x(\Z_\ell)=\Isom_k(\ul A(\ell))$ for all $\ell$.
  
  {\rm (2)} There is a natural bijection 
  \[ \Lambda_{x,N}(k)\simeq  
    G_{x}(\Q)\backslash G_{x}(\A_f)/U_N. \]
  
  {\rm (3)} $\mass[ \Lambda_{x.N}(k)]=\mass(G_{x}, U_N)$.
\end{thm}
\begin{proof}
  The statements (2) and (3) follow from the statement (1) and
  Theorems~\ref{22} and \ref{24}. We prove (1). Let $k_0$ be a
  finitely generated field for which $\ul A$ is defined and we may
  assume that $k$ is an algebraic closure of $k_0$. 
  
  If $\char k=0$, then $G_\ell:=\rho_\ell(\calG_{k_0})$ is in
  $Z(\GAut_{B}(V_\ell))$ for all 
  $\ell$. By Faltings' theorem \cite{faltings:end}, one has 
  \[ G_x(\Z_\ell)=\Aut_{k_0}(\ul A(\ell))=Z(\Aut_{k}(\ul
  A(\ell)), G_\ell). \]
  It follows that $G_x(\Z_\ell)=\Aut_{k}(\ul A(\ell))$. 
  This proves the case of \ch zero.

  If $\char k=p>0$, then replacing $\ul A$ by an isogeny we may assume
  that $k_0$ is a finite field, as $A$ is of CM-type (Lemma~\ref{311}
  and a theorem of Grothendieck \cite[p.~220]{mumford:av}). 
  Let ${\rm Frob}$ be the
  geometric Frobenius in $\calG_{k_0}$ and $\pi_A$ the relative Frobenius
  endomorphism on $A$. 
  By the $p$-adic version of Tate's theorem
  on endomorphisms over finite fields, it then suffices to show that
  $\pi_A$ lies in $Z(B)\otimes \Q_p$, the center of
  $\End^0_{B}(A[p^\infty])$. Since $A$ is of arithmetic type,
  one has $\pi_A=\rho_\ell({\rm Frob})\in Z(B)\otimes \Q_\ell$
  (Lemma~\ref{37}). Consider $\Q[\pi_A]$ and $Z(B)$ as linear
  subspaces of $\End^0(A)$; then $ \Q[\pi_A]=\Q_\ell[\pi_A]\cap \End^0(A)
  \subset Z(B)\otimes \Q_\ell\cap \End^0(A)=Z(B)$. Therefore, $\pi_A\in
  Z(B)$; this proves the statement (1). \qed  
\end{proof}


\section{Classification}
\label{sec:04}
In this section we classify abelian $B$-varieties of arithmetic type
up to isogeny. Due to Lemma~\ref{311}, it suffices to classify these objects
which are defined either over a number field or over a finite field. 
It is enough to consider the case where $B$ is a division algebra by
\ref{341}. From Lemma~\ref{41} to Proposition~\ref{43}, we assume that
$B$ is a division algebra. 


Let $\bbP$ be a prime field, $k$ be an algebraic closure of $\bbP$ and
$k_0$ be a finite extension of $\bbP$ in $k$. Let $\ul A$ be an
abelian $O_B$-variety over $k_0$. 

\begin{lemma}\label{41}
  If the positive involution $*$ on $B$ is of first kind, then $\ul A$ is of
  arithmetic type if and only if $\char k=p>0$ and $A$ is supersingular.  
\end{lemma}
\begin{proof}
  If $\ul A$ is of arithmetic type, then by Lemma~\ref{37} $G_\ell$ is
  contained in $\Q_\ell^\times$ after replacing $k_0$ by a finite
  extension. Then $\End^0(A)$ has dimension $4g^2$ by Tate's theorem. 
  This implies $\char k=p>0$ and $A$ is supersingular. The other 
  implication is obvious. \qed
\end{proof}

\subsection{} Lemma~\ref{41} classifies the abelian varieties of
  $B$-arithmetic type in the case of first kind.  
  Suppose that $*$ is of second kind. Let $K$ be the center
  of $B$ and $F$ be the maximal totally real subfield of $K$.  

\subsubsection{} Let $\iota_0:O_K\to \End(A)$ be the restriction of
$\iota$. Then $(A,\iota)$ is of $B$-arithmetic type if and only if
$(A,\iota_0)$ is of $K$-arithmetic type. Indeed, it follows from
Lemma~\ref{37} that the centers $Z(\GAut_{B}(V_\ell,\<\,,\>))$ and
$Z(\GAut_{K}(V_\ell,\<\,,\>))$ are the same. Therefore, the
classification is reduced to the case where $B$ is a CM field $K$. 

\subsubsection{}\label{422} The abelian variety $A$ is isogenous to 
$\prod_{j=1}^r A_j^{n_j}$, 
denoted by $A\sim \prod_{j=1}^r A_j^{n_j}$,
where each $A_j$ is a simple abelian variety and $A_i$ is not isogenous
to $A_j$ for $i\not = j$. If $\ul A$ is of $K$-arithmetic type, then we have
\[ \End^0_K(A)\otimes \Q_\ell\simeq \End_{K_\ell}(V_\ell). \]
Note that $V_\ell$ is a free $K_\ell$-module. The right hand side is
isomorphic 
to $M_n(K_\ell)$ and has $\Q_\ell$-dimension $n^2 d$, where $[K:\Q]=d$
and 
$2g=dn$. Put $B_j:=A_j^{n_j}$, $b_j=\dim B_j$ and let $m_j:=2b_j/d$. We
have 
\[\dim_\Q \End^0_K(B_j)\le \dim_{\Q_\ell} 
\End_{K_\ell}(V_\ell(B_j))=dm_j^2. \]
We also have $\End^0_K(A)=\prod_{j} \End_K^0(B_j)$.
From the dimensions of the abelian varieties $B_j$ and those of their
endomorphism algebras, we have
\[ \sum_{j=1}^r m_j=n,\quad n^2\le \sum_{j=1}^r m_j^2. \]
This shows $r=1$. We have shown that if $\ul A$ is of arithmetic type 
then it is isogenous to a self product of a simple factor. 

\begin{prop}\label{43}
  Suppose that $\char k=0$. Let $\ul A=(A,\iota)$ be an abelian
  $O_K$-variety with a CM field $K$. Then $\ul A$ is of $K$-arithmetic
  type if and only if
  $A\sim A_1^n$, where $A_1$ is a simple abelian variety with CM by a
  CM field $K_1$,
  and the image of the homomorphism $\iota: K\to \End^0(A)=M_n(K_1)$ 
  contains the center $K_1$.  
\end{prop}
\begin{proof}
  ($\Rightarrow$) Suppose that $A$ is of $K$-arithmetic type. The first 
  statement is proved in \ref{422}. We regard $K$ as a subfield of
  $\End^0(A)$ via $\iota$.
  Let $\wt K$ the composite of $K$ and 
  $K_1$. We want to show that $K=\wt K$. The centralizer of 
  $K$ in $M_n(K_1)$, same as that of $\wt
  K$, has $\Q$-dimension 
  \[ [\wt K:\Q] (\dim V / [\wt K:\Q])^2=4g^2/[\wt
  K:\Q].\] 
  While $\End_{K_\ell}(V_\ell)$ has $\Q_\ell$-dimension
  $4g^2/[K_\ell:\Q_\ell]$. It follows that $[\wt K:\Q]=[K:\Q]$ and 
  $K \supset K_1$. 

  ($\Leftarrow$) First, one has 
  $G_\ell\subset \End_{\End^0(A)}(V_\ell)$. Since
  $n[K_1:\Q]=\dim_{\Q_\ell} V_\ell$, the commutant  
  $\End_{\End^0(A)}(V_\ell)$ is $K_{1,\ell}$. The center of
  $\End_{K_1}(V_\ell)$ is $K_{1,\ell}$. So the Galois group $G_\ell$
  is contained in the center of $\End_{K_1}(V_\ell)$. 
  This shows that $A$ is of
  $K_1$-arithmetic type, particularly of $K$-arithmetic type. \qed
\end{proof}

Thereafter, the \ch of $k$ will
be $p>0$. 
We recall the definition of basic abelian varieties with additional
structures in the sense of Kottwitz (cf. \cite{kottwitz:isocrystals} and
\cite[p.~291, 6.25]{rapoport-zink}). 

\begin{defn} \label{44}
Let $W$ be the ring of Witt vectors over $k$ and $L$ be the fraction
field of $W$. Let $(B,*)$ remain as in \ref{21}.

{\rm (1)}  
  Let $(V_p,\psi_p)$ be a $\Q_p$-valued non-degenerate skew-Hermitian
  $B_p$-module, where $B_p:=B\otimes \Q_p$. A polarized abelian
  $O_B$-variety $\ul A$ over $k$ is said
  to be {\it related to $(V_p,\psi_p)$} if there is a 
  $B_p\otimes L$-linear isomorphism $\alpha:M(\ul A)\otimes_W L\simeq
  (V_p,\psi_p)\otimes L$ which preserves the pairings for 
  a suitable identification $L(1)\simeq L$, where $M(\ul A)$ is the
  covariant \dieu module with additional structures associated to $\ul A$. 

  Let $G':=\GAut_{B_p}(V_p,\psi_p)$ be the algebraic
  group over $\Q_p$ of $B_p$-linear similitudes. 
  A choice $\alpha$ gives rise to
  an element $b\in G'(L)$ so that one has an isomorphism of
  isocrystals with additional structures 
  $M(\ul A)\otimes L \simeq (V_p\otimes L, \psi_p, b({\rm id}\otimes
  \sigma))$.  
  The decomposition of $V_p\otimes L$ into
  isoclinic components induces a $\Q$-graded structure, and thus
  defines a (slope) homomorphism $\nu_{[b]}:\D\to G'$ over some finite
  extension $\Q_{p^s}$ of $\Q_p$, where $\D$ is the pro-torus over
  $\Q_p$ with character group $\Q$.    

{\rm (2)} A polarized abelian $O_B$-variety $\ul A$ is called {\it
  basic with respect to $(V_p,\psi_p)$} if 

(a) $\ul A$ is related to $(V_p,\psi_p)$, and

(b) the slope homomorphism $\nu$ is central. 

{\rm (3)}  $\ul A$ is called {\it basic} if it is basic
  with respect to $(V_p,\psi_p)$ for some skew-Hermitian space
  $(V_p,\psi_p)$.  
\end{defn}

\begin{lemma}\label{45} 
  Let $\ul A$  be a polarized abelian $O_B$-variety over $k$.
  The following statements are equivalent.

(a) $\ul A$ is basic.

(b) Let $Z$ be the center of $B$ and $Z_p=Z\otimes
\Q_p=\prod_{\bfp|p} Z_\bfp$ be the decomposition as a product of local
fields. Let $N=M(\ul A)\otimes_W L$ be the isocrystals with additional
structures associated to $\ul A$ and $N=\oplus_{\bfp|p} N_\bfp$ be the
decomposition with respect to the $Z_p$-action.  
Then each component $N_\bfp$ is isoclinic.   
\end{lemma}
\begin{proof}
  See a proof in 6.25 of \cite{rapoport-zink}.\qed
\end{proof}

One can use the statement (b) of Lemma~\ref{45} to check whether
an object $\ul A$ is basic. Note that the
statement (b) only depends on the underlying structure of $B$-action, 
not on the equipped polarization structure. 
This is also a property of those of arithmetic
type; see Lemma~\ref{38}. Indeed, we have

\begin{prop}\label{46} Let $(B,*)$ as in \ref{21}.
  An abelian $O_B$-variety $\ul A=(A,\iota)$ over $k$ is of arithmetic
  type if and only if it is basic. 
\end{prop}

\begin{proof}
  Write $A\sim \prod_j A_i^{n_i}$ into the decomposition up to isogeny
  with 
  respect to the decomposition $B=\oplus M_{n_i}(D_i)$ as in
  (\ref{34}). By (b) of Lemma~\ref{45} $A$ is basic if and only if
  each $A_i$ is basic. Therefore, we may assume that $B$ is a division
  algebra.  

  If $(B,*)$ is of first
  kind, then by Lemma~\ref{45} $\ul A$ is basic if and only if $\ul A$
  is supersingular. Then this follows from Lemma~\ref{41}. 

  Suppose that $(B,*)$ is of second kind. By Lemma 6.28 of
  Rapoport-Zink~\cite{rapoport-zink},
  $\ul A$ is basic if and only there is a finite field $k_0$ such that
  the relative Frobenius morphism $\pi_{A/k_0}$ lies in the center $K$
  of $B$. The latter statement is equivalent to that the Galois
  representation $\rho_\ell$ factors through the center
  $Z(\GAut_{B}(V_\ell,\<\,,\>))$; see the proof in
  Theorem~\ref{314}. This completes the proof. \qed
\end{proof}

By Proposition~\ref{46} and Theorem~\ref{314}, we get the following
result. This is the motivation of this work.

\begin{thm}\label{main}

  Let $(B,*)$ as in \ref{21}.
  Let $\ul A$ be a basic polarized abelian $O_B$-variety over an \ac
  field $k$ of \ch $p$, and let $N$ be any prime-to-$p$ positive
  integer. Then   

  {\rm (1)} There is a natural bijection 
  \[ \Lambda_{x,N}(k)\simeq  
    G_{x}(\Q)\backslash G_{x}(\A_f)/U_N. \]
  
  {\rm (2)} $\mass[ \Lambda_{x,N}(k)]=\mass(G_{x}, U_N)$.
\end{thm}

Consider the reduction $\calM\otimes \Fpbar$ modulo $p$ of a moduli
space associated to a PEL-datum $(B,*,V,\psi)$. Let
$x=(A,\lambda,\iota,\bar \eta)$ be a 
basic point in the moduli space $\calM\otimes \Fpbar$. 
By definition, any object $\ul A$
in $\calM(k)$ satisfies the condition $(I_\ell)$ of (i) for all
$\ell\neq p$. The set $\Lambda_{x,N}(k)$ essentially consists of
objects in $\calM(k)$ satisfying the condition $(I_p)$. The
condition ($Q$) is a mild technical condition which is subject to the
Hasse 
principle; see Remark below for details. 

\begin{remark}
(1) We come back to explain the meaning of the conditions (i), (ii)
    for the set $\Lambda_x(k)$ in Subsection~\ref{21} when $k=\Fpbar$. Let
    $x=\ul A_0=(A_0,\lambda_0,\iota_0)$ be a point fixed as
    before. Then the 
    space consisting of objects $\ul A$ that satisfy the condition (i)
    is the leaf $\calC(x)$ passing through $x$ in a moduli space
    $\calM\otimes \Fpbar$ of PEL-type. See \cite{oort:foliation} for
    the definition of leaves and detail discussions. The moduli space
    $\calM\otimes \Fpbar$ is the reduction modulo $p$ of a finite
    disjoint union of Shimura varieties of PEL-type due to the effect
    of the Hasse principle. When $x$ is basic, the dimension of the
    leaf $\calC(x)$ is zero. This is also, in fact, a necessary
    condition for 
    $\dim \calC=0$. 

(2) The condition (ii) is closely related to the prime-to-$p$ Hecke
    orbit of $x$. By definition, the prime-to-$p$ Hecke orbit
    $\calH^{(p)}(x)$ consists of objects $\ul A$ that satisfy the
    following conditions
\begin{itemize}
\item [${\rm (i)}'$] $\Isom_k(\ul
A_0(\ell),\ul A(\ell))\neq \emptyset$ for all $\ell\neq p$, and
\item [${\rm (ii)}'$] there is a prime-to-$p$ 
$O_B$-linear quasi-isogeny $\varphi: A_0 \to A$ which preserves the
polarizations. 
\end{itemize}
By definition we have inclusions $\calH^{(p)}(x)\subset
\Lambda_x(k)\subset \calC(x)$. When $x$ is basic, or more generally
$\ul A_0$ satisfies the condition
\begin{equation}
  \label{eq:47}
   \End_{O_B}(A_0)\otimes \Z_p\simeq \End_{O_B}(A_0[p^\infty]),
\end{equation}
then $\calH^{(p)}(x)= \Lambda_x(k)$. This follows from the weak
approximation of the reductive algebraic group $G_x$ (defined in
Subsection~\ref{21}). Abelian varieties with additional structures
that satisfy the condition (\ref{eq:47}) is called {\it
  $B$-hypersymmetric} 
in the sense of Chai. Hypersymmetric abelian varieties are studied in
Chai and Oort \cite{chai-oort:hyper}. 

(3) When $x$ is basic, $\Lambda_x(k)$ in general is strictly smaller
    than $\calC(x)$. That is, the condition (ii) is not automatically
    satisfied when the condition (i) is satisfied for an object $\ul
    A$. However, it is proved in Rapoport and Zink \cite[Chapter
    6]{rapoport-zink} that if $\ul A$ satisfies the condition (i),
    then there is an associated class $[\ul A]$ in $H^1(\Q,G_x)$ such
    that this class vanishes if and only if the condition (ii) is also
    satisfied. 
    Moreover, since $\ul A$ satisfies the condition (i), the class
    $[\ul A]$ 
    lies in the kernel $\ker^1(\Q,G_x)$ of the local-global map
\[ H^1(\Q,G_x)\to \prod_{v} H^1(\Q_v,G_x), \]
where $v$ runs through all places of $\Q$. Indeed, the map $\ul
A\mapsto [\ul A]$ gives a bijection between the set of quasi-isogeny
classes of objects $\ul A$ satisfying the condition (i) and the set
$\ker^1(\Q,G_x)$. In particular, when $G_x$
satisfies the Hasse principle, one has the quality
$\Lambda_x(k)=\calC(x)$. This is the case for the type C family of
Shimura varieties because $G_x$ is an inner form of the defining
reductive group $G_1$, which is semi-simple and simply connected. 
The reader may find that in \cite{yu:thesis}
and \cite{yu:mass_hb} the condition (ii) is not imposed; this is
because the Hasse principle for $G_x$ is satisfied. In general, we do
need to impose such a condition, which replaces the failure of the
Hasse principle. 

\end{remark}
  

\end{document}